\documentclass[10pt,a4paper]{amsart}
\usepackage{amsfonts,amsmath,amssymb}
\usepackage{epsfig}
\usepackage[super,nospace]{cite}
\let\refcite\citen

\newtheorem{dummy}{realdumb}

\newtheorem{theorem}[dummy]{Theorem}
\newtheorem{proposition}[dummy]{Proposition}

\newcommand{\bodymatter}{\maketitle}

\begin{document}

\title{3-MANIFOLDS AND 4-DIMENSIONAL SURGERY}

\author{M. YAMASAKI}

\address{Department of Applied Science, Okayama University of Science, 
Okayama, Okayama 700-0005, Japan, 
E-mail: masayuki@mdas.ous.ac.jp}

\begin{abstract}
Let $X$ be a compact connected orientable Haken 3-manifold with boundary,
and let $M(X)$ denote the 4-manifold $\partial(X\times D^2)$.
We show that if $(f,b):N\to M(X)$ is a degree 1 $TOP$ normal map
with trivial surgery obstruction in $L_4(\pi_1(M(X)))$,
then $(f,b)$ is $TOP$ normally bordant to a homotopy equivalence
$f':N'\to M(X)$.
Furthermore, for any $CW$-spine $B$ of $X$, we have a $UV^1$-map
$p:M(X)\to B$ and, for any $\epsilon>0$, 
$f'$ can be chosen to be a $p^{-1}(\epsilon)$-homotopy equivalence.
\end{abstract}


\bodymatter

\section{Introduction}
Hegenbarth and Repov\v{s} [\refcite{HR}] compared the controlled surgery exact
sequence of Pedersen-Quinn-Ranicki [\refcite{PQR}] with the ordianry surgery
sequence and observed the following:

\begin{theorem}[Hegenbarth-Repov\v{s}]
Let $M$ be a closed oriented $TOP$  4-manifold and $p:M\to B$
be a $UV^1$-map to a finite $CW$-complex such that the assembly map
\[
A:H_4(B;\mathbb{L}_\bullet)\to L_4(\pi_1(B))
\]
is injective.
Then the following holds:
if $(f,b):N\to M$ is a degree 1 $TOP$ normal map with trivial surgery obstruction
in $L_4(\pi_1(M))$, then $(f,b)$ is $TOP$ normally
bordant to a $p^{-1}(\epsilon)$-homotopy equivalence
$f':N'\to M$ for any $\epsilon>0$.
In particular $(f,b)$ is $TOP$ normally bordant to a homotopy equivalence.
\end{theorem}

\bigbreak
\noindent
{\bf Remarks.} 
(1) A map $f:N\to M$ is a {\it $p^{-1}(\epsilon)$-homotopy equivalence}
if there is a map $g:M\to N$ and homotopies $H:g\circ f\simeq 1_N$ and
$K:f\circ g\simeq 1_M$ such that all the arcs
\begin{align*}
&[0,1]\xrightarrow{H(x,-)} N \xrightarrow{~~~f~~~} M \xrightarrow{~~~p~~~} B\\
&[0,1]\xrightarrow{K(y,-)} M \xrightarrow{~~~p~~~} B
\end{align*}
have diameter $<\epsilon$.

\noindent
(2) $\mathbb{L}_\bullet$ is the 0-connective simply-connected surgery
spectrum [\refcite{Ra}].

\noindent
(3) The definition of $UV^1$-maps is given in the next section.
We have an isomorphism $\pi_1(M)\cong \pi_1(B)$.

\noindent
(4) This is true because the assembly map can be identified
with the forget-control map $F:H_4(B;\mathbb{L}_\bullet)\to L_4(\pi_1(M))$
which sends the controlled surgery obstruction to the ordinary
surgery obstruction.  By the injectivity of this map, the vanishing
of the ordinary surgery obstruction implies the vanishing of the controlled
surgery obstruction.

\bigbreak
For each torus knot $K$, 
Hegenbarth and Repov\v{s} [\refcite{HR}] constructed a 4-manifold $M(K)$
and a $UV^1$-map $p:M(K)\to B$ to a $CW$-spine $B$ of the exterior
of $K$ such that $A:H_4(B;\mathbb{L}_\bullet)\to L_4(\pi_1(B))$ is an isomorphism.
The aim of this paper is to extend their construction as follows.

Let $X$ be a compact connected orientable 3-manifold with nonempty boundary.
Then $M(X)=\partial(X\times D^2)$ is a closed orientable smooth 4-manifold
with the same fundamental group as $X$.
In fact, for any $CW$-spine $B$ of $X$, one can construct 
a $UV^1$-map $p:M(X)\to B$.

\begin{theorem}\label{main}
If $X$ is a compact connected orientable Haken 3-manifold with boundary,
and $B$ is any $CW$-spine of $X$, then there is a $UV^1$-map
$p:M(X)\to B$, and the assembly map $A:H_4(B;\mathbb{L}_\bullet)\to
L_4(\pi_1(B))$ is an isomorphism.
\end{theorem}

Thus we can apply Theorem 1 to these 4-manifolds.
Here is a list of such 3-manifolds $X$:
\begin{enumerate}
\item the exterior of a knot or a non-split link [\refcite{AFR}],
\item the exterior of an irreducible subcomplex of a triangulation of $S^3$
[\refcite{Ro1}].
\end{enumerate}

\bigbreak
The author recently learned that Qayum Khan proved the following [\refcite{Kh}].

\begin{theorem} [Khan]
Suppose $M$ is a closed connected orientable $PL$ 4-manifold
with fundamental group $\pi$ such that the assembly map
\[
A:H_4(\pi;\mathbb{L}_\bullet)\to L_4(\pi)
\]
is injective, or more generally, the 2-dimensional component
of its prime 2 localization
\[
\kappa_2:H_2(\pi;\mathbb{Z}_2)\to L_4(\pi)
\]
is injective.  Then any degree 1 normal map $(f,b):N\to M$
with vanishing surgery obstruction in $L_4(\pi)$ is normally
bordant to a homotopy equivalence $M\to M$.
\end{theorem}

In the examples constructed above, $X$'s are aspherical; so Khan's
theorem applies to the $M(X)$'s.

\medbreak
In \S\ref{UV}, we give a general method to construct $UV^m$-maps,
and finish the proof of Theorem \ref{main} in \S\ref{assembly}.

\section{Construction of $UV^{m-1}$-maps}
\label{UV}

A proper surjection $f:X\to Y$ is said to be $UV^{m-1}$ if,
for any $y\in Y$ and for any neighborhood $U$ of $f^{-1}(y)$ in $X$,
there exists a smaller neighborhood $V$ of $f^{-1}(y)$ such that
any map $K\to V$ from a complex of dimension $\le m-1$ to $V$ is
homotopic to a constant map as a map $K\to U$. 
A $UV^{m-1}$ map induces an isomorphism on $\pi_i$ for $0\le i<m$
and an epimorphism on $\pi_m$.
See [\refcite{La}, pp. 505--506] for the detail.

Let $X$ be a connected compact $n$-dimensional manifold
with nonempty boundary, and fix a positive integer $m$.
We assume that $X$ has a handlebody structure.
Recall from [\refcite{FQ}, p.136] that $X$ fails to have a handlebody
structure if and only if $X$ is an nonsmoothable 4-manifold.

Take the product $X\times D^m$ of $X$ with an $m$-dimensional disk $D^m$,
and consider its boundary $M(X)=\partial(X\times D^m)$,
which is an $(n+m-1)$-dimensional closed manifold.

Recall that a handlebody structure gives a $CW$-spine of $X$
[\refcite{KS}, p.107].
So, take any $CW$-spine $B$ of $X$:
there is a continuous map $q:\partial X\to B$ and $X$ is homeomorphic to
the mapping cylinder of $q$.  The mapping cylinder structure
extends $q$ to a strong deformation retraction $\overline{q}:X\to B$.
Define a continuous map $p:M(X)\to B$
to be the restriction of the composite map
\[
X\times D^m\xrightarrow{\text{projection}}X\xrightarrow{\overline{q}} B
\]
to the boundary.

\begin{proposition} For any $CW$-spine $B$ of $X$, 
$p:M(X)\to B$ is a $UV^{m-1}$-map.
\end{proposition}
\begin{proof}
First, let us set up some notations.
$M(X)$ decomposes into two compact manifolds with boundary:
\[
P=X\times S^{m-1}~,\qquad Q=\partial X\times D^m~.
\]
For any subset $S$ of $B$, define subsets $P_S\subset P$ and  $Q_S\subset Q$ by
\[
P_S=\overline{q}^{-1}(S)\times S^{m-1},\quad
Q_S=q^{-1}(S)\times D^m~.
\]
Then $p^{-1}(S)=P_S\cup Q_S$.

Let $b$ be a point of $B$ and take any open neighborhood $U$ of $p^{-1}(b)$
in $M(X)$.  
Since $M(X)$ is compact, the map $p$ is closed and hence 
there exists an open neighborhood $\widehat{U}$ of $b$ in $B$
such that $p^{-1}(\widehat{U})\subset U$.
Choose a smaller open neighborhood $\widehat{V}\subset \widehat{U}$ of $b$,
such that the inclusion map $\widehat{V} \to \widehat{U}$ is homotopic 
to the constant map to $b$, and set $V=p^{-1}(\widehat{V})$.

Suppose that $\varphi:K\to V$ is a continuous map from an 
$(m-1)$-dimensional complex.
We show that the composite map
\[
\varphi':K\xrightarrow{~~~\varphi~~~} V
\xrightarrow{\text{inclusion map}} U
\]
is homotopic to a constant map.

First of all, $Q_{\widehat{V}}$ has a core $q^{-1}(\widehat{V})\times\{0\}$
of codimension $m$, 
and, by transversality, we may assume that 
$\varphi:K\to P_{\widehat{V}}\cup Q_{\widehat{V}}$ misses
the core, and hence, we can homotop $\varphi$ to a map into $P_{\widehat{V}}$.
Since $P_{\widehat{V}}$ deforms into $\widehat{V}\times S^{m-1}$, 
we can further homotop $\varphi$ to a map into $\widehat{V}\times S^{m-1}$.
By the choice of $\widehat{V}$, $\varphi'$ is homotopic to a map into 
$\{b\}\times S^{m-1}$.
Pick any point $\overline{b}\in q^{-1}(b)$.
Then this map is homotopic to a map
\[
K\to \{\overline{b}\}\times S^{m-1}\subset \{\overline{b}\}\times D^m
\subset Q_{\widehat{U}}~.
\]
Therefore $\varphi'$ is homotopic to a constant map.
\end{proof}

\begin{proposition}
If $X$ has a handlebody structure, then 
$\pi_i(X)\cong\pi_i(M(X))$ for $i\le m-1$~,
and $\pi_m(X)$ is a quotient of $\pi_m(M(X))$.
\end{proposition}

\begin{proof}
This immediately follows from the proposition above,
but we will give an alternative proof here.

Take any handle decomposition of $X$:
\[
X = h_1\cup h_2\cup\dots\cup h_N~. 
\]
This defines the dual handle decomposition of $X$ on $\partial X$,
in which an $n$-handle of the original handlebody is a $0$-handle.
Since $X$ is connected, one can cancel all the $0$-handles
of the dual handle decomposition.
Thus we may assume that there are no $n$-handles in the
handlebody structure of $X$.

The handlebody structure of $X$ above
gives rise to a handlebody structure of $X\times D^m$:
\[
X\times D^m = h'_1\cup h'_2\cup\dots\cup h'_N~,
\]
where $h'_i=h_i\times D^m$ is a handle of the same index as $h_i$.
So there are only  $0$-handles up to  $(n-1)$-handles,
and the dual handle decomposition of $X\times D^m$ on $M(X)$
has no handles of index $\le m$.
The result follows.
\end{proof}

\bigbreak
\section{Proof of Theorem \ref{main}}
\label{assembly}

Roushon [\refcite{Ro2}] proved the following (among other things):

\begin{theorem} [Roushon]
Let $X$ be a compact connected orientable Haken 3-manifold.
Then the surgery structure set
$\mathcal{S}(X\times D^n \text{\it ~rel~} \partial)$ is trivial
for any $n\ge 2$.
\end{theorem}

The vanishing of $\mathcal{S}(X\times D^n \text{\it ~rel~} \partial)$
implies that the 4-periodic assembly maps [\refcite{Ra}]
\[
A:H_i(X;\mathbb{L}_\bullet(\mathbb{Z}))\to L_i(\pi_1(X))
\qquad\text{($i\in\mathbb{Z}$)}
\]
are all isomorphisms.
Since
\[
H_i(X;\mathbb{L}_\bullet(\mathbb{Z}))\cong H_i(X;\mathbb{L}_\bullet)
\]
for $i\ge\dim B$, the 0-connective assembly map
\[
A:H_4(X;\mathbb{L}_\bullet)\to L_4(\pi_1(X))
\]
is also an isomorphism.

Let $B$ be any $CW$-spine of $X$ and
let $p:M(X)\to B$ be the $UV^1$-map  constructed in the previous section.
Since $B$ is a deformation retract of $X$, the assembly map
\[
A:H_4(B;\mathbb{L}_\bullet)\to L_4(\pi_1(B))
\]
is an isomorphism.  This finishes the proof of Theorem \ref{main}.

\bigbreak
\section*{Acknowledgements}

This research was partially supported by 
Grant-in-Aid for Scientific Research 
from the Japan Society for the Promotion of Science.

I express my cordial thanks to Dusan Repov\v{s} for many helpful comments 
and to Jim Davis and Qayum Khan for their patient explanations of the results 
in [\refcite{Kh}] to me.

\end{document}